\documentclass[11pt]{amsart}
\usepackage{amsmath}
\usepackage{amssymb}
\usepackage{amsthm}
\usepackage{stmaryrd}
\usepackage[mathscr]{eucal}
\input xy
\xyoption{all}

\theoremstyle{plain}

\theoremstyle{definition}

\numberwithin{equation}{section}

\newcommand{\p}{\mathbb{P}}

\newcommand{\mc}{\mathcal}

\begin{document}

\title{Elliptic Curves with Bounded Ranks in Function Field Towers}
\author{Lisa Berger}
\address{Department of Mathematics\\
Stony Brook University\\
Stony Brook, NY 11794-3651}
\email{lbrgr@math.sunysb.edu}

\date{\today}


\maketitle

\section{Introduction}

We study the arithmetic structure of elliptic curves over $k(t)$, where $k$ is an algebraically closed field. In \cite{ShiodaAlgorithm} Shioda shows how one may determine rank of the N{\'e}ron-Severi group of a Delsarte surface--a surface that may be defined by four monomial terms. To this end, he describes an explicit method of computing the Lefschetz number of a Delsarte surface. He proves the universal bound of $56$ on the rank of an elliptic curve defined by an equation of the form $y^2 = x^3 + at^nx +bt^m$ over $k(t)$, where $k$ is an algebraically closed field of characteristic zero. In \cite{Shioda:Remarks} Shioda shows that the rank of $68$ is obtained for the curve $y^2 = x^3 +t^{360} +1$ over $\mathbb{C}(t)$. In recent work, Heinje \cite{Heijne:Delsarte} characterizes all Delsarte elliptic surfaces. He determines $42$ families of Delsarte elliptic curves and shows, through explicit computation, that $68$ is the maximal rank over $k(t)$, $k$ algebraically closed of characteristic zero. By relating a Delsarte surface to a Fermat surface, Shioda is able to exploit the relationship between divisor classes on his surface and the Mordell-Weil group of its generic fiber. In \cite{lisaTowers} the author describes a more flexible construction of elliptic surfaces. We explicitly construct families of surfaces, dominated by products of curves, with the additional property that they retain this DPC property under base extension. The N{\'e}ron Severi group of a product of curves may be expressed in terms of divisorial correspondences on the product, and Ulmer \cite{Ulmer:DPCT} utilizes this relationship to prove an explicit formula for the ranks of the Jacobians of the curves constructed in \cite{lisaTowers}. He produces elliptic curves with rank at least $13$ over $\mathbb{C}(t)$, and Occhipinti \cite{TommyThesis} produces an elliptic curve over $\bar{\mathbb{F}}_p(t)$ whose ranks over the fields $\bar{\mathbb{F}}_p(t^{1/d})$ grow at least linearly with $d$ prime to $p$. The goal of this note is to show that the large rank examples obtained via our construction are rare. We determine all elliptic curves obtained via the construction in \cite{lisaTowers}, and we find that, for all but finitely many families, the Mordell-Weil group of $E/k(t^{1/d})$ has rank zero, for each $d$ prime to the characteristic of $K=k(t)$, $k$ an algebraically closed field of arbitrary characteristic. \\

To state the main theorems, we first recall the construction and notation in \cite{lisaTowers}. Let $\mathcal{C}$ and $\mathcal{D}$ denote smooth, projective curves over a field $k$, and let $f$ and $g$ denote separable rational functions in $k(\mathcal{C})$ and $k(\mathcal{D})$, respectively. We have a canonically defined rational map: $\mathcal{C} \times _k \mathcal{D} \dashrightarrow \p^1_k$, $P \mapsto [f(P):g(P)],$ defined away from the locus of points $f=g=0$ and $f=g=\infty$. A blow-up of this locus resolves the map to a morphism from the often singular surface in $\mathcal{C} \times \mathcal{D} \times \p^1$, defined by the vanishing of $tf - g$, where $t = \frac{T}{S}$, $T$ and $S$ coordinates on $\p^1$. Let $\mathcal{S}$ denote a smooth, proper minimal model of this surface, with generic fiber $X_{f,g}$, a curve over $K = k(t)$. By construction, $\mathcal{S}$ is DPCT: it is dominated by a product of curves in towers of non-constant field extensions of the form $t \mapsto t^d$, $d$ prime to the characteristic of $k$. That the surface is DPC is clear; it is birational to $\mathcal{C} \times \mathcal{D}$. That this property is retained in towers is detailed in \cite{lisaTowers}. Let $m:= \text{deg}(f)$ and $n := \deg(g)$, $m_i$, $m'_{i'}$ the orders of the zeroes and poles of $f$, $n_j$, $n'_{j'}$ the orders of the zeroes and poles of $g$.

\subsection{Theorem}(\cite{lisaTowers}, \cite{Ulmer:DPCT})
\emph{Assume that the orders of zeros and poles of $f$ and $g$ have no common divisor and that they are relatively prime to the characteristic of $K = k(t)$. Then the generic fiber $X$ of a smooth projective model $S$ of the surface defined by the vanishing of $tf(x)-g(y)$ is an absolutely irreducible curve of geometric genus:}

\emph{$$g = mg_{D} + ng_{C} + (m-1)(n-1)- \sum_{(i,j)}\delta(m_i,n_j) - \sum _{(i',j')}\delta(m'_{i'},n'_{j'}),$$
}
\emph{where $g_D$ and $g_C$ denote the genera of the curves $\mathcal{D}$ and $\mathcal{C}$, respectively, and $\delta(a,b) = \frac{(a-1)(b-1)}{2}+\frac{((a,b)-1)}{2}$, and the sums are taken over all pairs $(i,j)$, $(i',j')$.} \\

 Let $\mathcal{C} = \mathcal{D} = \p^1$. Take rational functions $f$ and $g$ with $\text{div}(f) = \sum _{i=1}^k m_ia_i - \sum _{i'=1}^{k'} m' _{i'} a' _{i'}$ and $\text{div}(g) = \sum _{j=1}^{\ell} n_j b_j - \sum _{j' =1}^{\ell '} n' _{j'} b' _{j'}$, with all $a _i$, $b_j$, $a _{i'}$, and $a _{j'} \in k$, and with $a_i$, $a_{i'}$ all distinct and $b_{j}$, $b '_{j'}$ all distinct. Assume $(m,n) = 1$, and write $rm = \sum m_i = \sum m'_{i'}$ and $rn = \sum n_j = \sum n_{j'}$. Then the generic fiber, $X_{f,g}$, of the surface constructed above is a bidegree $(rm, rn)$ curve birational to the curve defined by the equation:
$tf(x)-g(y) = 0$. The main work we present in this note is an analysis of those partitions of $(rm, rn)$, the multiplicities of the zeros and poles of $f$ and $g$, for which our construction yields an absolutely irreducible curve with geometric genus one, and we obtain the following:

\subsection{Theorem}

\begin{enumerate}

\item

\emph{Let ${E}_{f,g}$ denote an elliptic curve over $k(t)$, constructed as above: the generic fiber of a smooth, proper model of the surface $tf-g \in \mathcal{C} \times \mathcal{D} \times \p^1$. Assume also that $rm := \deg(f) \leq \deg(g) = :rn$. Then, for all but finitely many bidegrees $(rm,rn)$, with (m,n)=1, $f$ has exactly one zero and one pole.}\\

\item

\emph{Let $K = k(t)$, $k=\bar{k}$, and let $E_{f,g}$ denote an elliptic curve over $K$, with defining equation as in the preceding statement: $\mc{C} = \mc{D} = \p^1$, $f$ has exactly one zero and one pole. Let $d$ range over non-negative integers, prime to the characteristic of $K$. Then the rank of the Mordell-Weil group of $E/k(t^{1/d})$ is zero.}\\

\end{enumerate}The proof of part one is computational and consists of an analysis of our genus formula, in the case of genus one. Along the way we give explicit models for the finitely many families of curves that are not of this form. Part two is a corollary to this classification theorem and to an explicit rank formula in \cite{Ulmer:DPCT}. \\

It is a pleasure to acknowledge the work of Erick Galinkin, a former Stony Brook undergraduate, who carried out some initial computations for this project. Thanks are also due to Tommy Occhipinti and Doug Ulmer for comments, suggestions and encouragement.

\section{Genus one partitions}

\subsection{} Take $\mathcal{C} = \mathcal{D} = \p^1$, construct the curve defined by $tf(x)-g(y)$ as above, and continue to assume in what follows that $m \leq n$. Set $\delta _0 := \sum _{i,j} \delta (m_i,n_j)$, $\delta _{\infty}:= \sum _{i',j'} \delta (m'_{i'}, n'_{j'})$, and $\delta := \delta _0 + \delta _{\infty}$. Our goal is to impose singularities with multiplicities to ensure that the smooth model $X_{f,g}$ has geometric genus one. We first explicitly determine the maximum obtainable value for $\delta_{0}$ and for $\delta_{\infty}$; we denote by $\delta _{\text{max}}$ this maximum value, and we show, without loss of generality, that a genus one curve may only be obtained when $\delta _0 = \delta _{\text{max}}$ or when $\delta _0 = \delta _{\text{max}}-\frac{r}{2}$. Finally, we describe the defining equations of all families of genus one curves obtained through our construction. Let $k$ and $k'$ denote the numbers of zeros and poles of $f$, $\ell$ and $\ell '$ the numbers of zeros and poles of $g$. \\

\subsection{Lemma}

\emph{Given positive integers $r$, $m$ and $n$ and a partition $([\{m_i\}],[\{n_j\}])$ of the bidegree $(rm,rn)$. The maximum possible value for $\delta _0$ is $\delta _{\text{max}} := \frac{r^2mn -rm -rn +r}{2}$.}

\begin{proof}

We have $$\delta _0 = \sum _{i=1}^k \sum _{j=1}^{\ell} \delta (m_i,n_j) = \sum _{i,j} \frac{(m_i-1)(n_j-1) + (m_i,n_j) - 1}{2},$$ and

$$\delta _0 = \frac{r^2mn - \ell rm - krn +\sum _{i,j} (m_i,n_j)}{2},$$

so, for fixed $r$, $m$ and $n$, we find the maximum possible value of $$D:=\sum _{i=1}^k \sum _{j=1}^{\ell} (m_i,n_j) - \ell rm - krn.$$\\

 When $\ell = k = 1$, we have $\sum _{i,j} (m_i,n_j) = (rm,rn)=r$, and $D=r - rm - rn$. We show that no larger value of $D$ may be obtained by increasing $\ell$ or $k$, the numbers of parts of our partitions. In what follows we suppose an increase in $k$. The argument is identical if we instead assume an increase in $\ell$.\\

 Re-ordering terms if needed, consider a partition: $m'_k + m''_k = m_k$ of $m_k$. We show that, for each $j$, $(m_k,n_j) + n_j \geq (m'_k,n_j) + (m'' _k,n_j)$. First suppose $n_j | m'_k$. If $n_j$ also divides $m''_k$ then we have equality. Otherwise, since $(m''_k, n_j)$ divides $m_k$, the inequality follows. If both $(m'_k,n_j)$ and $(m''_k,n_j) < n_j$ then their sum is bounded by $n_j$, and the strict inequality holds.\\

From this we obtain, $\sum _{j=1}^{\ell} (m_k, n_j) \geq \sum _{j=1}^{\ell} ((m'_k,n_j) + (m'' _k,n_j)) - rn$.\\

This yields

$$\sum _{i=1}^k\sum _{j=1}^{\ell} (m_i, n_j) \geq \sum _{i=1}^{k-1} \sum _{j=1}^{\ell}(m_k, n_j) + \sum _{j=1}^{\ell} ((m'_k,n_j) + (m'' _k,n_j)) - rn,$$

and

$$\sum _{i=1}^k\sum _{j=1}^{\ell} (m_i, n_j) - \ell rm - krn \geq \sum _{i=1}^{k-1} \sum _{j=1}^{\ell}(m_i, n_j) + \sum _{j=1}^{\ell} ((m'_k,n_j) + (m'' _k,n_j)) - \ell rm - (k+1)rn.$$

No larger value for $D$ may be obtained by increasing the number of elements in the partitions; the maximum value for $D$ is $r - rm -rn$, and the maximum value for $\delta _0$ and for $\delta _{\infty}$ is as claimed.\\

\end{proof}

\subsection{}

To obtain genus one we must choose partitions of $rm$ and $rn$ so that $\delta = 2\delta_{\text{max}} - r$. Indeed, letting $g_a$ denote the arithmetic genus, we have $g_a - 2\delta_{\text{max}} + r = (rm-1)(rn-1) - (r^2mn - rm - rn + r) + r = 1$. Assume without loss of generality that $\delta _{0} \geq \delta _{\infty}$. In the remainder of this section we find that a genus one partition is obtained only when $\delta _0 = \delta _{\infty} = \delta_{\text{max}}-\frac{r}{2}$ and when $\delta _0 = \delta_{\text{max}}$, $\delta _{\infty} = \delta_{\text{max}} - r$. We show that, for all but finitely many bidegrees $(rm,rn)$, we require $k = k' = 1$ to obtain genus one, and we determine all partitions that yield genus one.\\

We have
\begin{equation}
\delta _0 = \frac{r^2mn - \ell rm - krn -\sum _{i,j} (m_i,n_j)}{2}, \label{Eq:0}
\end{equation}
and if we assume $\delta _0 = \delta _{\text{max}} -\frac{r}{2}$ then we obtain the relation

\begin{equation}
(\ell - 1)rm + (k-1)rn = \sum (m_i, n_j) \leq \min \{\ell rm, krn\}. \label{Eq:1}
\end{equation}

We use the upper bound in \ref{Eq:1} to prove the following:

\subsection{Proposition}

\emph{Suppose $\delta _0 = \delta _{\infty} = \delta _{\text{max}} - \frac{r}{2}$. Then $(m,n) = (1,n)$ and $r=2$.}

\begin{proof}
\begin{itemize}
\item We assume first that $\ell$, $k \neq 1$. From the upper bound in \ref{Eq:1} we obtain $(\ell -1)rm \leq rn$ and $(k-1)rn \leq rm$. Combining these yields $(k-1)(\ell-1)rm \leq rm$, and this implies that $\ell = k = 2$. Making this substitution in \ref{Eq:1} we have

\begin{equation}
rm + rn = (m_1, n_1) + (m_1, n_2) + (m_2, n_1) + (m_2, n_2) \leq \min \{ 2rm, 2rn\}.\label{Eq:2}
\end{equation}

The upper bound in \ref{Eq:2} now implies that $m=n$ and, since $(m,n) = 1$, our bidegree is $(r,r)$. The equality in \ref{Eq:2} becomes: $2r = (m_1, n_1) + (m_2, n_1) + (m_1, n_2) + (m_2, n_2)$. Since, for $j=1, 2$, we have $\sum_i (m_i, n_j) \leq r$, each sum is exactly $r$. Hence, $(m_i, n_j) = m_i = n_j$, and all summands are equal. When each summand is $1$, so that the common divisor is one, we obtain an irreducible $(2,2)$ curve. Several families of $(2,2)$ curves are analyzed in \cite{lisaTowers}, \cite{TommyThesis} and \cite{Ulmer:DPCT}. Otherwise, for all $i$ and $j$, we have $(m_i,n_j) = \frac{r}{2} > 1$. Hence, in order to obtain an irreducible curve we now determine the complementary partitions $[\{m'_{i'}\}]$, $[\{n'_{j'}\}]$ of $(rm, rn) = (r,r)$ which yield $\delta _{\infty} = \delta _{\text{max}}-\frac{r}{2}$, satisfying $(m'_1, \cdots m' _{i'} , n'_1, \cdots n'_{j'},\frac{r}{2}) =1$. \\

From the upper bound in \ref{Eq:1}, assuming an $(r,r)$ curve, we find that the only possible partitions are of the form $\ell ' = k' = 2$ and $\ell ' = 2$, $k'=1$. (Since $m=n=1$, we need not consider the symmetric case $\ell ' = 1$, $k'=2$.) In the first case, as above, $\frac{r}{2} = (m'_1,m'_2,n'_1,n'_2,r)$. Hence, our bidegree is $(2,2)$. In the second case we obtain $r = (r, n'_1) + (r, n'_2)$. It follows, since $n'_1 + n'_2 = r$, that $(r, n'_j) = n'_j$, for each $j$. If $n'_1=1$ then $n'_2 = 1$, since $n'_2 | r$ and $r = 1 + n'_2$. So we have a $(2,2)$ curve. Otherwise, assume $(n'_1, n'_2) = 1$ but suppose, for some positive integers $k_1$ and $k_2$, that $r = n'_1k_1 = n'_2k_2 = n'_1 + n'_2$. If $r = n'_1n'_2 = n'_2+n'_2$, then $r=4$, and $2 = (n'_2,n'_2,r)$. Otherwise we must have $r > 4$, and we have $r > n'_1n'_2$, since each $n'_i$ divides $r$, and since $(n'_1,n'_2) = 1$. However, for $n'_1 + n'_2 = r$, $n'_i \neq 1$, we have $n'_1n'_2 > r$, a contradiction. So we obtain only bidegree $(2,2)$ curves when $\ell$, $k \neq 1$.\\

\item Assuming now that $k = 1$ for the first partition, again setting $\delta _0 = \delta _{\text{max}} - \frac{r}{2}$ yields:

$$(\ell - 1)rm = (rm, n_1) + (rm, n_2) + \cdots + (rm, n_{\ell}) \leq \min \{\ell rm, rn\}.$$

  The only possible set of summands is $rm + rm + \cdots + rm + \frac{rm}{2} + \frac{rm}{2}$. To ensure that the common divisor of the summands is one, we assume $rm = 2$. Since $r$ is even, $r=2$ and $m=1$, and we obtain families of $(2, 2n)$ curves.

\end{itemize}
\end{proof}

  We note that, except for the $(2,2)$ case described above, we have proved that, whenever $\delta _0 = \delta _{\infty} = \delta _{\text{max}} - \frac{r}{2}$ our genus one $(2,2n)$ models are determined by partitions of the form:\\

\begin{center}
\begin{tabular}{|l||c||c|} \hline

$(rm,rn)$     & & $[\{m_i\}]$ $[\{n_j\}]$, $[\{m'_{i'}\}]$ $[\{n'_{j'}\}]$\\ \hline

$(2,2n)$     & & $[2][2r_1, \cdots 2r_{\ell-2}, 2r_{\ell -1} +1, 2r_{\ell} +1], [2][2r'_1, \cdots 2r'_{\ell '-2}, 2r'_{\ell '-1} +1, 2r'_{\ell '} +1]$ \\ \hline

\end{tabular}\\
\end{center}
\vspace{.1in}

We next show that the only other way to obtain a genus one curve is by imposing singularities so that, without loss of generality, $\delta_{0} = \delta _{\text{max}}$.

\subsection{Proposition}

\emph{Suppose $a < \frac{r}{2}$ and let $\delta _0 = \delta_{\text{max}} - a$. Then $a=0$, $k=1$, and $(rm, n_1, \cdots , n_{\ell}, r) = r$}\\

\begin{proof}

Substituting $\delta _0 = \delta _{\text{max}} - a$ into equation \ref{Eq:0} gives:

\begin{equation}
(\ell-1)rm + (k-1)rn +r -2a = \sum (m_i,n_j) \leq \min \{\ell rm, krn\}. \label{Eq:3}
\end{equation}

We first note that either $\ell$ or $k$ must be equal to one: Since $a < \frac{r}{2}$, we have $r-2a >0$. Hence, if both $\ell$ and $k$ were greater than one, we would have $(\ell -1)rm + (k-1)rn +r-2a$ exceeding the upper bound in \ref{Eq:3}.\\

 Assuming $k=1$ and $\ell \geq 1$ in \ref{Eq:3}, we have $(\ell -1)rm + (r-2a) = (rm, n_1) + (rm, n_2) + \cdots (rm, n_{\ell})$. One possible solution is $(rm, n_i) = rm$, for $i = 1, \cdots (\ell - 1)$, and $(rm, n_{\ell}) = r - 2a$. With this solution $rm$ divides $n_i$, for $i = 1, \cdots, (\ell -1)$, and since $r$ divides $rn = \sum _{i=1}^{\ell} n_i$, it follows that $r$ divides $n_{\ell}$. Since $r$ also divides $rm$, $r | (r-2a)$, which is $(rm, n_{\ell})$. Since $r-2a$ is positive, it follows that $a=0$, and $\delta _0 = \delta_{\text{max}}$. Hence, $r$ divides each element of $\{rm, n_1, \cdots, n_{\ell}\}$\\

We also observe that there is no other set $\{(rm,n_j)\}$ satisfying $(rm, n_1) + \cdots (rm, n_{\ell}) = (\ell -1)rm + r$. Indeed, suppose for some $j$ that $(rm,n_j) < rm$.  We then have $(rm,n_j) \leq \frac{rm}{2}$. Hence, if two or more terms in our sum are each less than $rm$, we cannot sum to $(\ell -1) +r$. \\

Finally, since we assume $m \leq n$, the equality \ref{Eq:3} is not satisfied for $\ell =1$, $k \geq 1$.\\

\end{proof}

\subsection{}
We next determine the partitions of $(rm,rn)$ yielding $\delta _{\infty} = \delta _{\text{max}}^c := \delta _{\text{max}}-r$. Further, we are only interested in those partitions that satisfy $(m'_1, \cdots, m'_{k'}, n'_1, \cdots n'_{\ell '}, r) = 1$. Assuming $\delta _{\infty} = \delta _{\text{max}}^c$, (now writing $m_i$, $n_j$, $k$ and $\ell$ for $m'_{i'}$, $n'_{\ell '}$, $k'$ and $\ell '$), we have:

\begin{equation}
(\ell -1)rm + (k-1)rn -r = \sum(m_i, n_j) \leq \min \{\ell rm, krn\}. \label{Eq:4}
\end{equation}

Except for the case where $k=1$, there exist finitely many values of $\ell$ and $k$ that satisfy this relation. We will consider each of these cases and determine all corresponding bidegrees.  Toward this end, we have the following:

\subsection{Proposition}\label{Prop:lk}

\emph{Suppose $\delta _{\infty} = \delta _{\text{max}}^c$. Then:
}

\begin{enumerate}

\item \emph{$k=1$ and $\ell > 1$ or}
\item \emph{$k=2$ and $\ell = 2$, $3$, or $4$ or}
\item \emph{$k=\ell = 3$ or}
\item \emph{$\ell = 1$ and $k=2$ or $3$.}
\end{enumerate}

\begin{proof}

When $\ell = k = 1$ we have $\delta _{\infty} = \delta _{\text{max}}$, so $k=1$ implies $\ell > 1$ and $\ell = 1$ implies $k >1$. We next show that, for $k > 1$, $\ell$ is bounded above by $4$. From the upper bound in \ref{Eq:4} we obtain the relations $(k-1)n -1\leq m$ and $(\ell - 1) m -1 \leq n$. Combining these we obtain:

$$(\ell-1)m - 1 \leq \frac{m+1}{k-1},$$ and

$$\ell \leq \frac{k}{m(k-1)} + \frac{1}{k-1} + 1.$$

That $\ell \leq 4$ follows from the second inequality above, and this bound is obtained only when $k=2$ and $m=1$.\\

Beginning again with the bound in \ref{Eq:4}, we have

\begin{equation}
(\ell-1)(k-1)n - (\ell -1) \leq n + 1. \label{Eq:100}
\end{equation}

From this we obtain $(\ell k -\ell-k)n \leq \ell$, and we consider three cases.\\

 \emph{Case 1}: When $\ell k - \ell - k = 0$, we have $\ell k = \ell + k$, so $\ell = k = 2$.\\

 \emph{Case 2}: When $\ell k - \ell - k < 0$ we have $\ell (k-1) < k$. Either $k=1$ or $\ell < \frac{k}{k-1} <2$, so $\ell = 1$.\\

 \emph{Case 3}: Last, take $\ell k - \ell - k > 0$. Then, from \ref{Eq:100}, we obtain $n \leq \frac{\ell}{\ell k - \ell -k}$, so we determine those $\ell$ and $k$ for which $\frac{\ell}{\ell k - \ell -k} \geq 1$. Setting $\ell \geq \ell k - \ell -k$ we obtain $k \leq \frac{2\ell}{\ell -1} = 2 + \frac{2}{\ell -1}$. From this inequality it follows that $\ell = k = 3$ or $k \leq 2$.\\

It remains to show that $\ell = 1$ implies $k = 2$ or $k=3$. Substituting $\ell = 1$ into \ref{Eq:4} we have:

\begin{equation}
(k-1)rn - r= \sum (m_i,n_j) \leq \min \{rm, krn\}.\label{Eq:300}
\end{equation}

From the upper-bound in \ref{Eq:300} we have $(k-1)rn-r \leq rm$, so $n \leq \frac{m+1}{k-1}$. Since we also assume $m \leq n$, we have $m \leq n \leq \frac{m+1}{k-1}$. From this it follows that $k \leq 3$. When $k=3$ we find that $m=n=1$, so we obtain bidegree $(r,r)$. This case is identical to the case where $\ell = 3$ and $k=1$, so we will not consider this case below.

\end{proof}

We next examine each of the cases $(1)$-$(4)$ in Proposition \ref{Prop:lk}, and we determine all corresponding bidegrees. We show that, except for the cases where $k=1$, there are finitely many bidegrees satisfying $\delta _{\infty} = \delta _{\text{max}}^c$, under the additional assumption that $(m_i, \cdots, m_k, n_1, \cdots n_{\ell}, r) = 1$. We continue with our assumption that $m \leq n$.

\subsection{Proposition}\label{Prop:l=k=2}

\emph{Suppose $\delta _{\infty} = \delta _{\text{max}}^c$ that $\ell = k = 2$, and that $(m_1,m_2,n_1,n_2,r)=1$. Then the only possible bidegrees are: $(2,4)$, $(2,3)$, $(3,4)$, $(3,6)$, $(4,4)$, and $(4,6)$.}\\

\begin{proof}

When $\ell = k = 2$, the formula in \ref{Eq:4} becomes:

\begin{equation}
rm + rn - r = \sum_{i,j} (m_i, n_j) \leq \min{\{2rm, 2rn\}}. \label{Eq:5}
\end{equation}

The upper bound on the sum in \ref{Eq:5} limits the possible values for $m$ and $n$. Indeed, first set $rm + rn -r \leq 2rm$. Then $rn - r \leq rm$, so $n-1 \leq m$. Then, taking $rm + rn -r \leq 2rn$, we have $m-1 \leq n$. Combining these, taking $m \leq n$, we have $m=n$ or $m+1=n$.\\

\subsubsection{Case $1$: $n=m+1$} We first assume $n=m+1$ in \ref{Eq:5} and obtain $2rm = \sum (m_i, n_j)$. To attain this sum we must have, for each $j$, $(m_1, n_j) + (m_2, n_j) = rm$, since $rm$ is an upper bound on the sum of these two terms. Further, since, for each $i$ and $j$, $(m_i,n_j)$ is bounded by $m_i$, and since $m_1 + m_2 = rm$, we conclude for each $i$ and $j$ that $(m_i,n_j)$ is exactly $m_i$. This means that $m_i | n_j$ for each $i$ and $j$, so $m_i$ divides  $n_1 + n_2=rn$ for each $i$. We show next that each $m_i$ divides $n$.\\

Since $(m_i,r)$ divides $m_i$, $(m_i,r) | n_j$ for each $i$ and $j$, so $(m_i,r)$ divides each partition summand. Hence we must take $(m_i,r) = 1$. Combined with the fact that each $m_i$ divides $rn$, it follows that each $m_i$ divides $n=m+1$. When $m_1 = m_2$ we find that this common value divides $rm$ and $rn$, hence it divides $r$. From this we must assume either $m_1 \neq m_2$, or $m_1 = m_2 = 1$. If $m_1 = m_2=1$ then $rm=2$. When $m=2$ and $r=1$ we obtain bidegree $(2,3)$. If $m=1$ and $r=2$ we obtain bidegree $(2,4)$. \\

We assume now that $m_1 \neq m_2$, and since $m_1$ and $m_2$ divide $n=m+1$, we have $rm = m_1 + m_2 < 2m+2 \leq 2m+1$. From this it follows that $m=1$ and $r=3$, or $r=2$, or $r=1$.

\begin{itemize}

\item When $m=1$ and $r=3$ we obtain bidegree $(3,6)$. \\

\item When $r=2$ the bidegree is $(2m, 2m+2)$, and we have $m_1 + m_2 + 2 = 2n$. Since $m_1 | (m+1)$ and $m_2 | (m+1)$, and since we assume $m_1 \neq m_2$, we have, without loss of generality, $m_1 \leq (m+1)$ and $m_2 \leq \frac{m+1}{2}$. This yields $m \leq 3$. When $m=3$ we obtain bidegree $(6,8)$, and a partition with $\delta _{\infty} = \delta_{\text{max}}^c$ leaves a common divisor in the summands.  When $m=1$ we obtain bidegree $(2,4)$, and when $m=2$ we have bidegree $(4,6)$. \\

\item Finally assume $r=1$, so we have an $(m, m+1)$ curve, and formula \ref{Eq:5} becomes $2m = \sum(m_i,n_j)$. Since, for each $i$ and $j$, $(m_i,n_j) = m_i$, $m_i |(n_1+n_2)$. So $m_i |(m+1)$, which is equivalent to $m_i | (m_1 + m_2 + 1)$. From this we obtain $m_1|(m_2+1)$.  It follows that either $m_1 = m_2 = 1$ or that $m_1=1$ and $m_2 = 2$. Since we assume $m_1 \neq m_2$, we are in the latter case, and we obtain bidegree $(3,4)$.

\end{itemize}

\subsubsection{Case $2$: $m=n=1$} When $m=n=1$ we have an $(r,r)$ curve, and $\sum(m_i, n_j) = r$. When the two partitions of $r$ are identical, so that $m_1=n_1$ and $m_2 = n_2$, we have $\sum(m_i,n_j) = m_1 + m_2 + (m_1,n_2) + (m_2,n_1) > r$, a contradiction. So we assume distinct partitions $\{[m_i]\}$, $\{[n_j]\}$ of $r$, and we may also assume without loss of generality that $n_1 < m_1 \leq m_2 < n_2$. This gives us the inequalities $(m_1, n_1) \leq n_1$, $(m_2, n_1) \leq n_1$, $(m_1,n_2) \leq m_1$, and $(m_2,n_2) \leq m_2$. We set $d:=(m_2,n_2)$ and consider four cases.\\

\begin{itemize}

 \item Suppose first that $d \geq \frac{r}{3}$. Since $m_1 \leq m_2$, we have $m_2 \geq \frac{r}{2}$. Since $m_2$ is a multiple of $d$ satisfying $m_2 \geq \frac{r}{2}$, we take $m_2 = 2d  \geq  \frac{2r}{3}$. Then $n_2 \geq \frac{3r}{3} = r$, a contradiction, so $d < \frac{r}{3}$.\\

\item Suppose next that $\frac{r}{4} < d < \frac{r}{3}$. Since $m_2 \geq \frac{r}{2}$, $d \neq m_2$, and we assume $m_2 = 2d$ and that $n_2 = 3d$. Then $m_1 = r - 2d$ and $n_1 = r - 3d$. We have $\frac{3r}{4} < 3d$, so $n_1 < \frac{r}{4}$. It follows that $(n_1,m_1) < \frac{r}{4}$. Then $(n_1,m_2) < \frac{r}{8}$; this is because $(n_1, m_2) \neq (n_1, m_1)$, and because $(n_1, m_2) | n_1 < \frac{r}{4}$. Since $\frac{2r}{4} < 2d = m_2$, we have $m_1 < \frac{r}{2}$. Note also that $(m_1,n_2) \neq m_1$; otherwise we would have $(m_1,n_1)$ dividing $(m_1,m_2)$. So  $(m_1,n_1)$ would divide $m_1$, $n_1$ and $n_2$, and hence also $m_2$, since $m_1 + m_2 = n_1 + n_2 = r$. So we would have a common divisor in the partition. (If the common divisor is one, then $r=2$, and the $(2,2)$ curves have already been considered.) Since $(m_1,n_2) \neq m_1$, we have $(m_1,n_2) \leq \frac{r}{4}$. Then, since $(m_2,n_2) < \frac{r}{3}$, we have $\sum (m_i,n_j) < \frac{r}{8} + \frac{r}{4} + \frac{r}{4} + \frac{r}{3} = \frac{23r}{24} < r$, a contradiction. \\

\item In the case where $d = \frac{r}{4}$, we obtain $\delta_{\text{max}}^c$ via an $\ell = 2 = k$ partition of a $(4,4)$ curve with the partition $[2,2]$, $[3,1]$.\\

 \item Last, take $(m_2,n_2) < \frac{r}{4}$, and suppose that $(m_1, n_1) = n_1$ Then $(m_2,n_1)$ divides $m_1$, $m_2$, and $n_1$. Hence, it divides $r$ and also $n_2$, and there is a common divisor in the partition. We assume the common divisor is one, so $(m_2,n_1) = 1$, and each other term is strictly less than $\frac{r}{4}$. It follows that $r < 4$, and this is not possible, since we assume here that the partitions are distinct. Otherwise, for each $i$, we have $(m_i, n_1) < n_1 < \frac{r}{2}$, so each $(m_i,n_1)$ is less than $\frac{r}{4}$. Similarly, we have $(m_1,n_2)$ strictly less than $\frac{r}{4}$, and $\sum (m_i,n_j) < r$, a contradiction.\\

\end{itemize}

Tracing through this proof of Proposition $2.8$, we obtain the following genus one partitions:

\begin{center}
\begin{tabular}{|l||c|} \hline
$(rm,rn)$               & $[\{m_i\}]$ $[\{n_j\}]$, $[\{m'_{i'}\}]$ $[\{n'_{j'}\}]$\\ \hline
$(2,3)$                &[2][3], [1,1][2,1]\\
$(2,3)$                &[2][2,1], [1,1][2,1]\\ \hline
$(2.4)$                &[2][4], [1,1][2,2]\\
$(2.4)$                &[2][4], [1,1][3,1]\\
$(2.4)$                &[2][2,2], [1,1][2,2]\\
$(2.4)$                &[2][2,2], [1,1][3,1]\\\hline
$(3,6)$                &[3][6], [2,1][2,4]\\
$(3,6)$                &[3][3,3], [2,1][2,4]\\\hline
$(4,6)$                &[4][6], [3,1][3,3]\\
$(4,6)$                &[4][4,2], [3,1][3,3]\\\hline
$(3,4)$                &[3][4], [2,1][2,2]\\
$(3,4)$                &[3][3,1], [2,1][2,2]\\\hline
$(4,4)$                &[4][4], [3,1][2,2]\\\hline

\end{tabular}
\end{center}

\end{proof}

\subsection{Proposition}

\emph{Suppose $\delta _{\infty} = \delta _{\text{max}}^c$, $\ell = k = 3$ and that $(m_1, m_2, m_3, n_1, n_2, n_3, r) = 1$. Then the bidegree is $(3,3)$.}

\begin{proof}

 Substituting $\ell = k = 3$ into \ref{Eq:4} we have

 \begin{equation}
 2rm + 2rn - r = \sum_{i,j} (m_i, n_j) \leq \min \{3rm, 3rn\}.\label{Eq:6}
 \end{equation}

From the upper bound in \ref{Eq:6} we obtain the inequalities $2rn - r \leq rm$ and $2rm -r \leq rn$. Combining these we have $n=m=1$, so we obtain an $(r,r)$ model and determine that $3r = \sum(m_i,n_j)$. We have, for each $i$, $\sum _j(m_i,n_j) \leq r$ and, for each $j$, $\sum _i (m_i,n_j) \leq r$. Hence, since the total sum of terms is $3r$, we have equalities: $\sum _j(m_i,n_j) = r$ and $\sum _i (m_i,n_j) = r$. For any fixed $j$, consider the sum $\sum _i(m_i, n_j)$. Each term is bounded by $m_i$, and $\sum _i m_i = r = \sum _i(m_i, n_j)$. Hence, each term $(m_i, n_j)$ is exactly $m_i$. So, for each $i$, $(m_i, n_j) = m_i$. Analogously, for fixed $i$, for each $j$, $(m_i, n_j) = n_j$. From this it follows that $(m_i, n_j)=m_i = n_j$ for all $i$ and $j$. Since we assume $(m_1, m_2, m_3, n_1, n_2, n_3, r) = 1$, we obtain $(m_i, n_j) = m_i = n_j = 1$. The bidegree is $(3,3)$ and we have the partition:

\begin{center}
\begin{tabular}{|l||c|} \hline
$(rm,rn)$               & $[\{m_i\}]$ $[\{n_j\}]$, $[\{m'_{i'}\}]$ $[\{n'_{j'}\}]$\\ \hline

$(3,3)$                &[3][3], [1,1,1][1,1,1]\\\hline

\end{tabular}
\end{center}

 \end{proof}

\subsection{Proposition}

\emph{Suppose $\delta _{\infty} = \delta_{\text{max}}^c$, $\ell = 3$ and $k=2$ and that $(m_1, m_2, m_3, n_1, n_2, r)=1$. Then the bidegrees are: $(2,3)$, $(2,4)$, $(3,6)$, and $(4,4)$.}

\begin{proof}

Substituting $\ell = 3$ and $k=2$ into \ref{Eq:4} we have

\begin{equation}
2rm +rn -r = \sum(m_i,n_j) \leq \min \{3rm, 2rn\} \label{Eq:7}
\end{equation}

From the upper bound in \ref{Eq:7} we obtain the inequalities $rn-r \leq rm$ and $2rm-r \leq rn$. Combining these we find that $m \leq 2$ and $n \leq 3$. In particular, $m=2$ and $n=3$ or $m=1$ and $n=1$ or $m=1$ and $n=2$. We analyze each of these cases.

\subsubsection{Case $1$: $m=2$, $n=3$.} Assuming $m=2$ and $n=3$ in \ref{Eq:7} we obtain the relation $6r = \sum (m_i,n_j)$. Combining this with the fact that, for each $i$, $\sum_j (m_i, n_j) \leq 3r$, we have the equality $\sum_j (m_i, n_j) = 3r$. Analogously, we have, for each $j$, $\sum _i(m_i, n_j) = 2r$. It follows that, for each $i$ and $j$, $(m_i,n_j) = r$. Indeed, suppose without loss of generality that the term $(m_1, n_1)$ is less than $r$. Then $(m_2, n_1) = 2r - (m_1,n_1) > r$. Since our summands are positive integers, $(m_2,n_1)| m_2$ implies $m_2 > r$. Since $m_1 + m_2 = 2r$, $m_1 < r$. Then $(m_1,n_2), (m_1,n_3) < r$, contradicting $\sum _j (m_1,n_j) = 3r$. Hence each summand is exactly $r$, so the common divisor of our partition is $r$. Since we assume $(m_1, m_2, m_3, n_1, n_2, r)=1$, we find that $r=1$ and the bidegree is $(2,3)$.\\

\subsubsection{Case $2$: $m=n=1$.} Substituting $m=n=1$ into \ref{Eq:7} we obtain bidegree $(r,r)$, and $\sum (m_i, n_j) = 2r$. It follows, for each $i$ and $j$, that $(m_i, n_j) = n_j$. Since $n_j$ divides $m_1$ and $m_2$, $n_j$ divides $m_1 + m_2 = r$. Hence, for each $j$, $n_j$ is an integer of the form $\frac{r}{a}$, and we determine those $a$, $b$, $c \in \mathbb{Z}^+$ satisfying $\frac{r}{a} + \frac{r}{b} + \frac{r}{c} = r$. The only solutions are $(a,b,c) = (2, 3, 6)$, $(3,3,3)$ and $(2, 4, 4)$.  Each of the first two triples corresponds to bidegree $(6,6)$. From the first we cannot find a corresponding partition of $r = m_1 + m_2$. From the second, the common divisor of the partition summands is $2$. From the last triple we obtain bidegree $(4,4)$.\\

\subsubsection{Case $3$: $m=1$, $n=2$.} Assuming $m=1$ and $n=2$ in \ref{Eq:7} we obtain $3r = \sum (m_i, n_j)$. Hence, for each $j$,  $\sum_i (m_i, n_j) = r$. It follows that each $(m_i, n_j) = m_i$. This is because $m_1 + m_2 = rm = r$ and because no divisor can exceed $m_i$. Hence, for each $i$ and $j$, $m_i | n_j$, and it follows that each $m_i$ divides $\sum n_j = 2r$. So, we have divisor sums of the form $\frac{2r}{a} + \frac{2r}{b} = r$, where $a$ and $b$ are non-negative integers. The only solutions are $a=3$, $b=6$ and $a=b=4$. The former corresponds to bidegree $(3,6)$ and the latter to bidegree $(2,4)$.\\

The bidegrees and partitions are summarized in the table:

\begin{center}
\begin{tabular}{|l||c|} \hline
$(rm,rn)$               & $[\{m_i\}]$ $[\{n_j\}]$, $[\{m'_{i'}\}]$ $[\{n'_{j'}\}]$\\ \hline

$(2,3)$                &[2][3], [1,1][1,1,1]\\
$(2,3)$                &[2][2,1], [1,1][1,1,1]\\\hline
$(4,4)$                &[4][4], [2,2][2,1,1]\\\hline
$(2,4)$                &[2][4], [1,1][2,1,1]\\
$(2,4)$                &[2][2,2], [1,1][2,1,1]\\\hline
$(3,6)$                &[3][6], [2,1][2,2,2]\\
$(3,6)$                &[3][3,3], [2,1][2,2,2]\\\hline

\end{tabular}
\end{center}

\end{proof}

\subsection{Proposition}

\emph{Suppose $\delta = \delta_{\text{max}}^c$, $\ell = 4$ and $k = 2$, and that $(m_1, m_2, n_1, n_2, n_3, n_4, r) = 1$. Then the bidegree is $(2,4)$.}

\begin{proof}

 Substituting $k=2$ and $\ell = 4$ into \ref{Eq:4} we have

 \begin{equation}
 3rm + rn - r = \sum(m_i,n_j) \leq \min \{4rm, 2rn\}. \label{Eq:8}
 \end{equation}

 From the upper bound in \ref{Eq:8} we obtain $rn - r \leq rm$ and $3rm -r \leq rn$. Combining these inequalities we find that $n=2$ or $n=1$. When $n=1$ we have $m < \frac{2}{3}$, so there is no corresponding bidegree. Setting $n=2$ in \ref{Eq:8} we obtain $m=1$, and the only possible bidegree has the form $(r, 2r)$. Making this substitution in \ref{Eq:8} yields $4r = \sum (m_in_j)$. Hence, for each $j$, $(m_1, n_j) + (m_2, n_j) = r$. From this it follows that $(m_i,n_j) = m_i$ for each $i$ and $j$. So, since $m_i | n_j$, it follows that $m_i | \sum n_j$, so $m_i | 2r$. We determine positive integers $a$ and $b$ with $\frac{2r}{a} + \frac{2r}{b} = r$. The only solutions $(a,b)$, up to reordering, are $(3,6)$ and $(4,4)$.  Reasoning as in the preceding section, one obtains the bidegrees $(2,4)$ and $(3,6)$. Only $(2,4)$ yields genus one. The partitions are:\\

\begin{center}
\begin{tabular}{|l||c|} \hline
$(rm,rn)$               & $[\{m_i\}]$ $[\{n_j\}]$, $[\{m'_{i'}\}]$ $[\{n'_{j'}\}]$\\ \hline

$(2,4)$                &[2][4], [1,1][1,1,1,1]\\
$(2,4)$                &[2][2,2], [1,1][1,1,1,1]\\\hline

\end{tabular}
\end{center}

\end{proof}

We next determine those bidegrees corresponding to $\delta _{\infty} = \delta _{\text{max}}^c$ when $\ell =1$ and  $k=2$. Still assuming $m \leq n$, we prove:  \\

\subsection{Proposition}
\emph{Suppose $\delta _{\infty} = \delta _{\text{max}}^c$, $\ell = 1$, $k = 2$ and $(m_1, m_2, rn, r) =1$. Then, the bidegrees are: $(2,3)$, $(2,4)$, $(3,4)$, $(3,6)$, $(4,6)$, and $(5,6)$.}

\begin{proof}

When $k=2$ the relation \ref{Eq:300} is

\begin{equation}
rn - r = (m_1, rn) + (m_2, rn) \leq \min \{rm, 2rn\}, \label{Eq:310}
\end{equation}

 \noindent
 so $n \leq m + 1$. Since the case $m=n$ has already been considered above, and since we assume $m \leq n$, we take $n=m+1$. Making this substitution in \ref{Eq:310} yields $rm = (m_1, rn) + (m_2, rn)$, so $m_1$ and $m_2$ each divide $rn = r(m+1)$. As in the proofs of the preceding two propositions, we determine positive integer solutions to: $\frac{m+1}{a} + \frac{m+1}{b} = m$, and we obtain the triples $(m,a,b)$: $(1, 3, 6)$, $(1, 4, 4)$, $(2, 2, 6)$, $(2, 3, 3)$, $(3, 2, 4)$ and $(5, 2, 3)$. \\

The first solution yields bidegree $(2,4)$ and the second bidegree $(3,6)$. The third integer triple corresponds to bidegree $(2,3)$ and the fourth to bidegree $(4,6)$. The fifth triple corresponds to bidegree $(3,4)$, and the last triple corresponds to bidegree $(5,6)$. We have the following partitions:\\

\begin{center}
\begin{tabular}{|l||c|} \hline
$(rm,rn)$               & $[\{m_i\}]$ $[\{n_j\}]$, $[\{m'_{i'}\}]$ $[\{n'_{j'}\}]$\\ \hline

$(2,3)$                &[2][3], [1,1][3]\\
$(2,3)$                &[2][2,1], [1,1][3]\\\hline
$(2,4)$                &[2][4], [1,1][4]\\
$(2,4)$                &[2][2,2], [1,1][4]\\\hline
$(3,4)$                &[3][4], [1,2][4]\\
$(3,4)$                &[3][3,1], [1,2][4]\\\hline
$(3,6)$                &[3][6], [1,2][6]\\
$(3,6)$                &[3][3,3], [1,2][6]\\\hline
$(4,6)$                &[4][6], [1,3][6]\\
$(4,6)$                &[4][4,2], [1,3][6]\\\hline
$(5,6)$                &[5][6], [2,3][6]\\
$(5,6)$                &[5][5,1], [2,3][6]\\\hline

\end{tabular}
\end{center}

\end{proof}

\subsection{} In the discussion above we determined all possible genus one bidegrees and partitions for which $k=1$, $k'>1$. We next determine all those obtained by setting $k=k'=1$. Still writing $\ell$ and $k$ in place of $\ell '$ and $k'$, substituting $k=1$ into \ref{Eq:4}, we have:

\begin{equation}
(\ell - 1) rm - r = (rm, n_1) + \cdots + (rm, n_{\ell}). \label{Eq:400}
\end{equation}

For simplicity, we first consider the case where $r=1$, and prove the following:

\subsection{Proposition}

\emph{Suppose $\delta _{\infty} = \delta _{\text{max}}^c$, $k=1$ and $r=1$. Then the bidegrees are $(2,n)$, $(3,n)$, $(4,n)$ and $(6,n)$, where $n$ may be any integer satisfying $n \geq m$ and $(m,n)=1$.}

\begin{proof}

  Assuming $r = 1$ in \ref{Eq:400} we obtain $(\ell - 1)m - 1 = (m, n_1) + \cdots (m, n_{\ell})$, and we first note that all but three of the terms $(m, n_j)$ must be equal to $m$. Indeed, supposing there are four terms less than $m$, we have $(\ell - 4)m + 4(\frac{m}{2}) \geq (\ell - 1)m - 1$, so $m \leq 1$, and our curve $X_{f,g}$ would not have genus one. Hence, the partition is of the form:

$$(\ell - 1)m - 1 = m + \cdots + m + (m, n_{\ell -2}) + (m, n_{\ell - 1}) + (m, n_{\ell}),$$

so we determine restrictions on the last three terms, and we need

\begin{equation}
2m - 1 = (m, n_{\ell -2}) + (m, n_{\ell - 1}) + (m, n_{\ell}).\label{Eq:420}
\end{equation}

We first assume that each term on the right hand side in \ref{Eq:420} is less than $m$. From this we obtain $\frac{3m}{2} \geq 2m-1$, which implies that $m \leq 2$. When $m=2$ we obtain $(2,n)$ curves, $2 \nmid n$. Since we assume $r=1$, we cannot have $m=1$, since in this case $X_{f,g}$ would be a rational curve.\\

We next assume the first term is $m$, leaving $m-1 = (m, n_{\ell - 1}) + (m, n_{\ell})$. One possible sum is $\frac{m}{2} + (\frac{m}{2} - 1)$, so we determine conditions for which $(\frac{m}{2} - 1) | m$. We have implicitly assumed $m$ is even, and the only solutions to the divisibility condition are $m=4$, and $m=6$. When $m=4$ we obtain $(4, n)$ curves, $n$ odd. When $m=6$ we obtain $(6, n)$ curves, $n \equiv 1 \text{ or } 5 \pmod 6$. Another possible sum is $\frac{m}{3} + (\frac{2m}{3} -1)$. In this case we have $m \equiv 0 \pmod 3$ and $(\frac{2m}{3} -1) | m$. These conditions imply that $m=3$ or that $m=6$. The case $m=6$ is identical to, (with partition symmetric to), the $m=6$ case above. When $m=3$ we obtain $(3,n)$ curves, $n \geq 2$. Setting the sum as $\frac{m}{4} + (\frac{3m}{4} -1)$ we obtain $m=4$, and this case has also been completed. There are no other partitions $\frac{m}{t} + (\frac{(t-1)m}{t} - 1)$ of $m-1$.\\

\end{proof}

\subsection{} In the table below we summarize the genus one partitions for $(rm,rn)$, curves, under the assumption that $r=1$ and $k=k'=1$:

\begin{center}
\begin{tabular}{|l||c||c|} \hline
$(rm,rn)$     &                      & $[\{m_i\}]$ $[\{n_j\}]$, $[\{m'_{i'}\}]$ $[\{n'_{j'}\}]$ \\ \hline
$(2,n)$      & $n \equiv 1 \pmod 2$ &$[2][2r_1, \cdots 2r_{{\ell}-1}, 2r_{\ell} + 1], [2][2r'_1, \cdots 2r'_{{\ell}' -3}, 2r'_{{\ell}'-2} + 1, 2r'_{{\ell}'-1}+1, 2r' _{{\ell}'}+1]$\\\hline
$(3,n)$      & $n \equiv 1 \pmod 3$ &$[3][3r_1, 3r_2, \cdots 3r_{{\ell}-1}, 3r_{\ell} + 1], [3][3r'_1, 3r'_2, \cdots 3r'_{{\ell}'-2}, 3r'_{{\ell}'-1}+2, 3r' _{{\ell}'}+2]$\\
             & $n \equiv 2 \pmod 3$ &$[3][3r_1, 3r_2, \cdots 3r_{{\ell}-1}, 3r_{\ell} + 2], [3][3r'_1, 3r'_2, \cdots 3r'_{{\ell}'-2}, 3r'_{{\ell}'-1}+1, 3r' _{{\ell}'}+1]$\\ \hline
$(4,n)$      & $n \equiv 1 \pmod 4$ &$[4][4r_1, 4r_2, \cdots, 4r_{{\ell}-1}, 4r_{\ell} +1], [4][4r'_1, 4r'_2, \cdots, 4r'_{{\ell}'-1}+2, 4r'_{{\ell}'} +3]$\\
             &  $n \equiv 3 \pmod 4$ &$[4][4r_1, 4r_2, \cdots, 4r_{{\ell}-1}, 4r_{\ell} +3], [4][4r'_1, 4r'_2, \cdots, 4r'_{{\ell}'-1}+2, 4r'_{\ell '} +1]$ \\ \hline
$(6,n)$      & $n \equiv 1 \pmod 6$ &$[6][6r_1, 6r_2, \cdots 6r_{{\ell}-1}, 6r_{\ell} + 1], [6][6r'_1, 6r'_2, \cdots 6r'_{\ell '-2}, 6r'_{\ell '-1}+3, 6r' _{\ell '}+4]$  \\
             & $n \equiv 5 \pmod 6$ &$[6][6r_1, 6r_2, \cdots 6r_{\ell-1}, 6r_{\ell} + 5], [6][6r'_1, 6r'_2, \cdots 6r'_{\ell '-2}, 6r'_{\ell '-1}+3, 6r' _{\ell '}+2]$  \\ \hline
\end{tabular}
\end{center}

\hspace{.1in}\\

We next assume $r > 1$, and we have:

\subsection{Proposition}

\emph{Suppose $\delta _{\infty} = \delta _{\text{max}}^c$, $k=1$, $r > 1$ and $(rm, n_1, \cdots n_{\ell}, r) = 1$. Then the bidegree is of the form: $(2,2n)$, $(3,3n)$, $(4,4n)$, $(4, 4s+2)$, $(6,6n)$, $(6, 6s+2)$, $(6, 6s+4)$ or $(6, 6s+3)$, where $s$ is a positive integer}.

\begin{proof}

Substituting $k=1$ into \ref{Eq:4} gives $(\ell -1)rm -r = \sum _{j=1}^{\ell}(rm,n_j)$, and we first note that we may not have more than four terms in the sum less than $rm$. This would give
$(\ell - 5)rm + \frac{5rm}{2} \geq (\ell - 1)rm - r$, and we would have $m \leq \frac{2}{3}$.

\begin{itemize}

\item In the case where we have exactly $4$ terms less than $rm$, we have $3rm - r \leq 2rm$, which implies $m \leq 1$. Setting $m=1$ in \ref{Eq:4} we have $2r = \sum _{j=1}^4 (r,n_j)$. Since we assume $(rm, n_1, \cdots n_{\ell}, r) = 1$, it follows that $r = 2$. We obtain a family of $(2, 2n)$ models, $n$ odd, with the partitions:

\begin{center}
\begin{tabular}{|l||c|} \hline
$(2,2n)$     & $[2][2r_1, \cdots 2r_{{\ell}-1}, 2r_{\ell}], [2][2r'_1, \cdots, 2r'_{\ell ' - 1}, 2r'_{\ell '-3}+1, 2r'_{\ell '-2} +1, 2r'_{\ell '-1}+1, 2r'_{\ell '} +1]$ \\ \hline

\end{tabular}
\end{center}

\hspace{.1in}\\

\item When exactly three terms are not equal to $rm$ it is sufficient to consider partitions that satisfy:

$$2rm - r = (rm, n_1) + (rm, n_2) + (rm, n_3).$$

 Since each term is less than $rm$ we have $\frac{3rm}{2} \geq 2rm - r$. Then $3m \geq 4m - 2$, so $m \leq 2$.

Setting $m=2$ we have:

$$3r = (2r, n_1) + (2r, n_2) + (2r, n_3),$$

and the only solution,comes from a $(2,n)$ curve, $r = 1$, since $(rm, n_1, \cdots, n_{\ell}, r) = r$. This case has been completed above. Setting $m=1$ we consider partitions that satisfy:

$$r = (r, n_1) + (r, n_2) + (r, n_3),$$

The only possible partitions are $(\frac{r}{3} + \frac{r}{3} + \frac{r}{3})$, $(\frac{r}{2} + \frac{r}{4} + \frac{r}{4})$ and $(\frac{r}{2} + \frac{r}{3} + \frac{r}{6})$; the only possible bidegrees are $(3,3n)$, $(4,4n)$, and $(6,6n)$, with the partitions: \\

\begin{center}
\begin{tabular}{|l||c||c|} \hline
$(3,3n)$     & & $[3][3r_1, \cdots, 3r_{\ell}], [3][3r'_1, \cdots 3r'_{\ell -3}, 3r'_{\ell -2}+1, 3r'_{\ell -1}+1, 3r'_{\ell}+1]$
\\ \hline

$(4, 4n)$    & & $[4][4r_1, \cdots, 4r_{\ell}], [4][4r'_1, \cdots 4r'_{\ell - 3}, 4r'_{\ell -2}+2, 2r'_{\ell -1}+1, 2r'_{\ell}+1]$\\ \hline

$(6, 6n)$    & & $[6][6r_1, \cdots, 6r_{\ell}], [6][6r'_1, \cdots 6r'_{\ell - 3}, 6r'_{\ell -2}+3, 6r'_{\ell -1}+2, 6r'_{\ell}+1]$ \\ \hline

\end{tabular}
\end{center}

\item Finally, take the case where exactly two terms are less than $rm$, so we consider partitions that satisfy:

$$rm - r = (rm, n_{\ell -1}) + (rm, n_{\ell}),$$

Reasoning as in the proof of the preceding proposition, we obtain, for the partition $\frac{rm}{2} + (\frac{rm}{2} - r)$, the restriction that $m=3$, $m=4$, or $m=6$. The first possibility, $m=3$, is only possible in the case where $r$ is even. In fact, we obtain new partitions only for $m=3$; when $m=4$ and $m=6$, we have $r=1$, which has been considered above. Further, when $m=3$ we have the restriction $r=2$, and we have bidegrees $(6, 6s+2)$ and $(6, 6s+4)$. \\

For the partition $\frac{rm}{3} + (\frac{2rm}{3} - r)$, we obtain the restriction $m=2$, $m=3$, or $m=6$, and we may have $m=2$ only if $r \equiv 0 \pmod 3$. We obtain new partitions only when $m=2$ and $r=3$, so we have the bidegrees $(6, 6s+3)$.\\

For the partition $\frac{rm}{4} + (\frac{3rm}{4} - r)$ we have $m=2$ or $m=4$, and $m=2$ is possible only if $r$ is even. Further, we obtain new partitions only when $m=2$ and $r=2$, and we have the bidegrees $(4, 4s+2)$. The partitions are:\\

\begin{center}
\begin{tabular}{|l||c|} \hline
$(6, 6s+2)$& $[6][6r_1, \cdots 6r_{\ell -1}, 6r_{\ell} +2], [6][6r'_1, \cdots 6r'_{\ell '-2}, 6r'_{\ell -1} +3, 6r'_{\ell}+5]$\\

 $(6, 6s+4)$& $[6][6r_1, \cdots 6r_{\ell -1}, 6r_{\ell} +4], [6][6r'_1, \cdots 6r'_{\ell '-2}, 6r'_{\ell -1} +3, 6r'_{\ell}+1]$ \\

  $(6, 6s+3)$ & $[6][6r_1, \cdots 6r_{\ell -1}, 6r_{\ell} +3], [6][6r'_1, \cdots 6r'_{\ell '-2}, 6r'_{\ell -1} +2, 6r'_{\ell}+1]$ \\ \hline

$(4, 4n+2)$ & $[4][4r_1, \cdots 4r_{\ell -1}, 4r_{\ell} +2], [4][4r'_1, \cdots 4r'_{\ell '-2}, 2r'_{\ell -1} +1, 2r'_{\ell}+1]$ \\ \hline
\end{tabular}
\end{center}

\end{itemize}

\end{proof}

\subsection{}Summarizing the main results above, we note that genus one partitions are obtained only from bidegrees $(2, N)$, $(3,N)$, $(4,N)$, $(5,6)$ and $(6,N)$. Further, all but finitely many of families come from those partitions that satisfy $k=k'=1$. The table below describes these exceptional families, those for which the defining function $f(x)$ does not have a unique zero and a unique pole.

\begin{center}
\begin{tabular}{|l||c||c||c|} \hline
& Bidegree                & Families & \\ \hline
1& $(2,2)$       &$tx(x-1)(y+1)(y-b)=(x+1)(x-a)y(y-1)$& $a$, $b \neq 0$, $1$\\ \hline

2& $(2,4)$       &$tx^2(y-1)(y-a)(y-b)(y-c) = y^2 (y-d)^2 (x-1)(x+1)$& $a$, $b$, $c \neq 0$,\\
          &               &                                                   &$d \neq 1$, $a$, $b$, $c$\\ \hline

3 & $(2,3)$      &$tx^2 (y-1)(y-a)(y-b) = y^2(y-d) (x-1)(x+1)$ & $a$, $b \neq 0$ \\
           &              &                                             &  $d \neq a$, $b$, $1$\\ \hline

4&$(3,3)$         &$tx^3(y-1)(y+1)(y-a) = y^3(x-1)(x+1)(x - b)$& $a$, $b \neq 0$, $1$, $-1$\\ \hline

5&$(3,4)$         &$tx^3 (y-1)^2 (y-a)^2 = y^3(y-b) (x-1)(x+1)^2$& $a \neq 0$, $b$ \\
          &                &                                              & $b \neq 1$ \\ \hline

6& $(3,6)$        &$t(y-a)^2(y-b)^2(y-1)^2 x^3 = y^3(y-d)^3 (x-1)^2(x+1)$& $a$, $b \neq 0$,\\
          &                &                                                      & $d \neq a$, $b$, $1$\\ \hline

7&$(4,4)$         &$t(y-1)^2(y+1)(y-a)x^4 = y^4 (x-1)^2(x+1)^2$&$a \neq -1$, $0$ \\ \hline

8&$(4,6)$        &$tx^4 (y-1)^3(y-a)^3= y^4(y-b)^2 (x-1)^3(x+1)$& $a \neq 0$, $b$ \\
 &                        &                                              &  $b \neq 1$\\ \hline

9&$(5,6)$                 &$tx^5 (y-1)^6 = (x-1)^3(x+1)^2 y^5 (y-a)$& $a \neq 1$\\ \hline
\end{tabular}
\end{center}

Note that, for each bidegree, a displayed family may admit degeneration and hence correspond to more than one of the partitions we determined in this section.\\

\subsection{Other Products $\mathcal{C} \times \mathcal{D}$}

Above we restrict to the case considered in \cite{lisaTowers}, setting $\mathcal{C} = \mathcal{D} = \p ^1$. Here we show that there are no other curves $\mathcal{C}$, $\mathcal{D}$ for which our construction yields a genus one curve at the base of the tower. When $\mathcal{C}$ and $\mathcal{D}$ are elliptic curves the genus is:

$$g = rm + rn + (rm-1)(rn-1) - \sum _{i,j} \delta(m_i, n_j) - \sum _{i',j'} \delta(m_i ', n_j ').$$

Setting $g = 1$, simplifying, we obtain:

$$(\ell + \ell ')rm + (k + k')rn = \sum _{i,j} (m_i, n_j) + \sum _{i',j'} (m_i ', n_j '),$$

where $\ell$, $\ell '$, $k$ and $k'$ are defined as above. We have already noted that the sum $\sum (m_i, n_j) + \sum (m_i ', n'_j) \leq  \min \{krn, \ell rm\}+ \min \{k' rn, \ell ' rm\}$. It follows that we cannot obtain a genus one model via our construction in this case, and an analogous argument shows that we cannot consider curves $\mathcal{C}$, $\mathcal{D}$ of higher genus. The only other case where a genus one curve could be obtained at the base of our construction would be for $\mathcal{C} = \p ^1$ and $\mathcal{D} = E$, an elliptic curve. In that case we obtain the restriction:

\begin{equation}
(\ell + \ell ')rm + (k + k' -2)rn = \sum _{i,j}(m_i,n_j) + \sum _{i', j'} (m_i ',n_j ').\label{Eq:500}
\end{equation}

From the upper bounds on each of the sums on the right hand side of equation \ref{Eq:500} one shows that $k = k' = 1$, and we have:

\begin{equation}
\ell rm = (rm, n_i) + \cdots + (rm,n_{\ell}) \text{  and  }\ell ' rm = (rm, n'_{i'}) + \cdots + (rm,n'_{\ell '}).
\end{equation}

 Each summand is $rm$, and this is possible only when $rm = 1$; otherwise we would not have an absolutely irreducible generic fiber. But when $rm=1$ we have a rational curve. We have proved, and now restate, Theorem $1.2$ $(1)$:

\subsection{Theorem}

\emph{Let $E_{f,g}$ denote the elliptic curve over $k(t)$, constructed as above, the generic fiber of a smooth, proper model of the surface $tf-g \subseteq \mathcal{C} \times \mathcal{D} \times \p ^1$. Assume $\deg(f) = rm \leq rn = \deg(g)$. Then, for all but finitely many bidegrees $(rm,rn)$, with $(m,n)=1$, $f$ has exactly one zero and one pole.} \\

\section{Bounded Ranks}

\subsection{}

We have shown that, for all but finitely many bidegrees, to obtain genus one curves via our construction we require $f(x)$ defined with exactly one zero and one pole. Assume as stated in the introduction that $k$ is an algebraically closed field. In this section we study our $(rm,rn)$ genus one curves over the fields $K=k(t)$. The main theorem is: except for the exceptional families in the table above, all of our elliptic curves $E/k(t^{1/d})$ have Mordell-Weil groups with rank zero, for all $d$ prime to the characteristic of $K$.

\subsection{} For a global field $K$, by the Mordell-Weil theorem, the group $E(K)$ is a finitely generated abelian group. One may also consider curves over the fields $k(t)$ where $k$ is an arbitrary field, and the group of $k(t)$ points of the Jacobian variety $J_K:=J(X_{f,g})$ need not be finitely generated. Let $A$ denote an abelian variety over the field $K$. One defines the $K/k$-trace of $A$ to be an abelian variety $B/k$ with a $K$-homomorphism $\tau: B \otimes _k K \rightarrow A$, satisfying the following universal property: If $C/k$ is another abelian variety with homomorphism $\psi: C _k \otimes _k K \rightarrow A$, then we have a homomorphism $\tau ':C \otimes _k K \rightarrow B \otimes _k K$, and the following commutative diagram:

$$\xymatrix{
C_k \otimes _k K \ar[d]_{\tau '} \ar[dr]^{\psi}&\\
B_k \otimes _k K \ar[r]_{\tau} & A_K.}$$

That is, $B$ is the largest abelian variety, defined over $k$, with $B \times _k K$ mapping to $A$ as above. In this work we are interested in the case where $A_K$ is the Jacobian variety $J_K$, as above. We define the Mordell-Weil group $MW(J_K):=J(K)/\tau (B(k))$, and by the Lang-N{\'e}ron Theorem \cite{LangNeron}, this quotient is a finitely generated abelian group. See \cite{ConradTrace} for a complete discussion of the $K/k$ trace of an abelian variety over $K$ and the Lang-N{\'e}ron Theorem.

\subsection{} Let $\mathcal{S}_d$ denote the base change $t \mapsto t^d$ of the surface $\mathcal{S}$ described in the introduction. In a recent preprint Ulmer shows that the surface $\mathcal{S}_d$ is a birational model of the quotient $(\mathcal{C}_d \times \mathcal{D}_d)/\mu _d$, where the curves $\mathcal{C}_d$, $\mathcal{D} _d$ are smooth, projective models of: $w^d = f(x)$ and $v^d = g(y)$, and where $\mu _d$ acts via: $(x, w) \mapsto (x, \zeta_d w)$, $(y, v) \mapsto (y, \zeta_d v)$, \cite{Ulmer:DPCT}. He considers a birational model $\mc{X}_d$ of $S_d$, and a smooth, proper morphism $\pi _d: \mc{X}_d \rightarrow \p^1$, which factors through $(\mathcal{C}_d \times \mathcal{D}_d)/\mu _d$, and he uses the geometry of this construction to determine an explicit formula for the ranks of the Mordell-Weil groups, as defined above, of the Jacobians of our curves over $k(t)$. Key to the rank formula is our construction of an elliptic surface as a birational model of a product variety, so that the N{\'e}ron-Severi group of the surface may be expressed in terms of divisorial correspondences on $\mathcal{C} \times \mathcal{D}$. The rank formula follows from this, combined with the Shioda-Tate formula and a thorough analysis of the geometry in our construction. In this section we use Ulmer's rank formula to bound the ranks of the elliptic curves described in the preceding section. Combined with the classification theorem, we find that the large rank examples in \cite{lisaTowers}, \cite{TommyThesis} and \cite{Ulmer:DPCT} are rare: there are finitely many bidegrees $(rm,rn)$ for which our construction yields elliptic curves with non-zero rank over the fields $k(t^{(1/d)})$ when $k$ is algebraically closed.

\subsection{} We briefly discuss Ulmer's rank formula and refer the reader to \cite{Ulmer:DPCT} for further details.  Let $f_{d,v}$ denote the number of irreducible components in the fiber of $\pi _d:\mathcal{X}_d \rightarrow \p ^1$, over the closed point $v$. Define $c_1(d) := \sum _{v \neq 0, \infty} (f_{d,v} -1).$ When $k$ is algebraically closed this becomes

 $$c_1(d) = d\sum_{v\neq 0, \infty}(f_{1,v}-1).$$

 Let $P_i$ and $P'_{i'}$ denote the zeros and poles of $f$, $Q_j$, $Q'_{j'}$ the zeros and poles of $g$, and let $t_{d,i,j}$ and $t'_{d,i',j'}$ denote the numbers of closed points of the surface $(C_d \times D_d)/ \mu_d$ over the points $(P_i, Q_j)$ and $(P_i', Q_j')$, respectively. Let $f'_{d,0}$ and $f'_{d,\infty}$ denote the number of irreducible components in the fibers of $(C_d \times D_d)/ \mu_d \dashrightarrow \p^1$ lying over $0$ and $\infty$, respectively.  \\

   Define $$c_2(d) := \sum _{i,j} t_{d, i, j} + \sum _{i,j} t' _{d, i', j'} - f'_{d,0} - f'_{d,\infty} + 2.$$\\

Note also that, in our construction, the covers $\mathcal{C}_d$ and $\mathcal{D}_d$ are often reducible. In the case where the base curves $\mathcal{C}$ and $\mathcal{D}$ are both rational, we let $e_{d,f}$ and $e_{d,g}$ denote the number of irreducible components of $C_d$ and $D_d$, respectively.  We write $\mathcal{C}_d'$ and $\mathcal{D}'_d$ for the smooth, proper models of $w^{d/{e_{d,f}}} = (f(x))^{1/e_{d,f}}$ and $v^{d/e_{d,g}}=(g(y))^{1/{e_{d,g}}}$, respectively. When we take for our constant field an algebraic closure of $k$, we have $t_{d,i,j} = \gcd(m_i, n_j, d)$ and $t'_{d, i', j'} = \gcd(m_i', n_j', d)$. Our formula becomes:

$$c_2(d) := \sum _{i,j} \gcd(m_i, n_j, d) + \sum _{i,j} \gcd(m'_{i'}, n'_{j'}, d) - \sum _i (m_i,e_{d,g}) - \sum_j (n_j,e_{d,f}) - \sum_{i'}(m'_{i'},e_{d,g}) - \sum_{j'} (n'_{j'}, e_{d,g}) + 2.$$\\

The constant $c_2(d)$ varies with $d$; it is clearly periodic, hence bounded. We have the following:

 \subsection{Theorem} (\cite{Ulmer:DPCT}). \emph{Let $\mathcal{C}$ and $\mathcal{D}$ denote smooth projective curves over and algebraically closed field $k$, $f \in k(\mathcal{C})$, and $g \in k(\mathcal{D})$, separable rational functions. Let $X_{f,g}$ denote a smooth model of the curve $tf -g$, constructed as above. Write $J_{\mc{C}_d}$ and $J_{\mc{D}_d}$ for the Jacobians of the curves $\mc{C}_d$ and $\mc{D}_d$, and write $J=\text{Jac}(X_{f,g})$ and $J_d = J/k(t^{1/d})$.  With notation as above, the rank over $k(t^{1/d})$ of the Mordell-Weil group of the Jacobian of $X_{f,g}$ is:}

\emph{$$\text{Rank} MW(J_d) = \text{RankHom}_{k-av}(J_{\mathcal{C}'_d}, J_{\mathcal{D}'_d})^{\mu_{d/({{e_{d,f}} \cdot e_{d,g}})}} - c_1(d) + c_2(d),$$}

where $\text{Hom}_{k-av}(J_{\mathcal{C}'_d}, J_{\mathcal{D}'_d})^{\mu_{d/({e_{d,f}} \cdot e_{d,g})}}$ denotes those homomorphisms commuting with the action of the group $\mu _{d/({e_{d,f} \cdot e_{d,g}})}$.\\

In the remainder of this section we show that, for ``most" elliptic curves arising in our construction, the $\text{Rank Hom}_{ab-v}(J_{\mathcal{C}}, J_{\mathcal{D}})^{\mu _d}$ and $c_1(d)$ terms in Ulmer's formula are both zero, and that the ranks are bounded in towers of function field extensions. To complete the proof of Theorem $1.2$ we first focus on the rank of $\text{Hom}_{k-av}(J_{\mathcal{C}'_d}, J_{\mathcal{D}'_d})^{\mu_{d/({e_{d,f}} \cdot e_{d,g})}}$, writing $\mu '_d$ for $\mu_{d/({e_{d,f}} \cdot e_{d,g})}$

 \subsection{Lemma}

 \emph{With notation as in the statement of the theorem, let ${X}_{f,g}$ denote the curve over $k(t)$, constructed as above: the generic fiber of a smooth, proper model of the surface $tf-g \in \mathcal{C} \times \mathcal{D} \times \p^1$, and assume also that $f$ has exactly one zero and one pole. Then $\text{RankHom}_{\text{ab-v}}(J_{C'_d}, J_{D'_d})^{\mu '_{d}}=0$}, and the invariant $c_1(d) = 0$ for all $d$.

\begin{proof}

In this case the curve $C'_d$ is a smooth, projective model of: $w^{d/e_{d,f}} = x^{rm/e_{d,f}}$, a rational curve. The Jacobian $J_{\mathcal{C'}_d}$ is trivial; $\text{RankHom}_{\text{ab-v}}(J_{C'_d}, J_{D'_d})^{\mu '_d} = 0$; and there is no contribution from this term to the rank of $E_{f,g}(k(t^{1/d}))$.\\

To prove the second part of the Lemma, we suppose first that $d=1$.  Since $\text{Rank Hom}_{ab-v}(J_{\mathcal{C}}, J_{\mathcal{D}})^{\mu _d} = 0$. The rank formula reduces to $\text{Rank} (X_{f,g}(K)) = -c_1(1) + c_2(1)$. But $c_2(1) = (\ell - 1)(k-1) + (\ell ' -1)(k' - 1)$, and we assume, $k = k' = 1$, so $c_2(1) = 0$. Since $\text{Rank} (X_{f,g}(K)) \geq 0$, $c_1(1) = 0$, and since $c_1(d)$ is linear in $d$, the Lemma follows.
\end{proof}

Hence, the Mordell-Weil rank of the Jacobian of our curve, in the case where $f$ has exactly one zero and pole,  is determined by the function $c_2(d)$. To bound this rank over the fields $K_d:=k(t^{1/d})$, completing the proof of Theorem $1.2 (2)$, we prove the following:

\subsection{Lemma}

 \emph{Let ${E}_{f,g}$ denote the elliptic curve over $k(t)$, constructed as above: the generic fiber of a smooth, proper model of the surface $tf-g \in \mathcal{C} \times \mathcal{D} \times \p^1$, $\mc{C} = \mc{D} = \p^1$, where $f$ has exactly one zero and one pole. Then the invariant $c_2(d) = 0$ for all $d$.}

\begin{proof}

Recall that we consider the curves over $k(t^{1/d})$, $k=\bar{k}$. In this case, the invariant $c_2(d)$ is:

\begin{align*}
c_2(d) = & \sum _{i,j} (m_i,n_j,d) + \sum_{i', j'} (m_{i'}, n_{j'}, d) - (\sum_j (d, n_j, rm) + (d, n_1, \cdots , n_j, n'_1, \cdots n '_{j'}, rm)) \\
& - (\sum_{j'} (d, n'_{j'}, rm) + (d, n'_{1}, \cdots , n'_{j'}, n_1, \cdots n _{j}, rm)) + 2,
\end{align*}

and this simplifies to:

$$c_2(d) =  \sum _{i,j} (m_i,n_j,d) + \sum_{i', j'} (m_{i'}, n_{j'}, d) - \sum_j (d, n_j, rm) - \sum_{j'} (d, n'_{j'}, rm).$$

\noindent

But $\sum _{i,j} (m_i,n_j,d) = \sum _j (rm, n_j, d)$, and $\sum_{i', j'} (m'_{i'}, n'_{j'}, d) = \sum _{j'} (rm, n'_{j'}, d)$,

so the invariant $c_2(d) = 0$.

\end{proof}

This proves the second part of Theorem $1.2$, which we restate here:

\subsection{Theorem}

\emph{Let $K = k(t)$, $k$ an algebraically closed field, and let $E_{f,g}$ denote an elliptic curve over $K$, the generic fiber of the surface $tf-g \in \mc{C} \times \mc{D} \times \p^1$, and assume that $f$ has exactly one zero and one pole. Let $d$ range over non-negative integers, prime to the characteristic of $K$. Then the rank of the Mordell-Weil group of $E/k(t^{1/d})$ is zero.}

\section{Remarks}

\subsection{}We show that, for all but finitely many bidegrees, all elliptic curves arising via our construction have rank zero in the towers $E(k(t^{1/d}))$. It is clear that the combinatorial argument could be extended to classify our Jacobians of higher dimension.

\subsection{}In \cite{AvanziZannier}, Avanzi and Zannier give a complete classification of genus one curves defined by equations of the form $f(x)=g(y)$, $f(x)$, $g(x) \in K[x]$, where $K$ is a field of characteristic zero, under the assumption that $\gcd (\deg{f}, \deg{g}) = 1$. In our classification of genus one curves we repeat some of the results of Avanzi-Zannier, but for our construction we are able to say more. First, we consider rational functions $f(x)$ and $g(x)$. Second, we have a stronger irreducibility result for our curves, allowing us to remove the assumption that $\deg{f}$ and $\deg{g}$ are relatively prime.\\

\bibliography{Dissertation}

\begin{thebibliography}{Ulm11}

\bibitem[AZ01]{AvanziZannier}
Roberto~M. Avanzi and Umberto~M. Zannier.
\newblock Genus one curves defined by separated variable polynomials and a
  polynomial {P}ell equation.
\newblock {\em Acta Arith.}, 99(3):227--256, 2001.

\bibitem[Ber08]{lisaTowers}
Lisa Berger.
\newblock Towers of surfaces dominated by products of curves and elliptic
  curves of large rank over function fields.
\newblock {\em J. Number Theory}, 128:3013--3031, 2008.

\bibitem[Con06]{ConradTrace}
Brian Conrad.
\newblock Chow's {$K/k$}-image and {$K/k$}-trace, and the {L}ang-{N}\'eron
  theorem.
\newblock {\em Enseign. Math. (2)}, 52(1-2):37--108, 2006.

\bibitem[Hei11]{Heijne:Delsarte}
Bas Heinje.
\newblock The maximal rank of elliptic {D}elsarte surfaces.
\newblock {\em http://arxiv.org/abs/1011.2340v2}, 2011.

\bibitem[LN59]{LangNeron}
S.~Lang and A.~N{\'e}ron.
\newblock Rational points of abelian varieties over function fields.
\newblock {\em Amer. J. Math.}, 81:95--118, 1959.

\bibitem[Occ10]{TommyThesis}
Thomas Occhipinti.
\newblock {\em {M}ordell-{W}eil groups of large rank in towers}.
\newblock PhD thesis, University of Arizona, 2010.

\bibitem[Shi86]{ShiodaAlgorithm}
Tetsuji Shioda.
\newblock An explicit algorithm for computing the {P}icard number of certain
  algebraic surfaces.
\newblock {\em American Journal of Mathematics}, 108(2):415--432, 1986.

\bibitem[Shi92]{Shioda:Remarks}
Tetsuji Shioda.
\newblock Some remarks on elliptic curves over function fields.
\newblock {\em Ast\'erisque}, (209):12, 99--114, 1992.
\newblock Journ{\'e}es Arithm{\'e}tiques, 1991 (Geneva).

\bibitem[Ulm11]{Ulmer:DPCT}
Douglas Ulmer.
\newblock On {M}ordell-{W}eil groups of {J}acobians over function fields.
\newblock {\em http://arxiv.org/abs/1002.3310v3}, 2011.

\end{thebibliography}
\bibliographystyle{alpha}

\end{document}